\documentclass{elsarticle} 

\usepackage{natbib}       

\usepackage{graphicx,hyperref}
\usepackage{booktabs} 
\usepackage{algorithm}
\usepackage{algpseudocode}

\usepackage[T1]{fontenc}
\usepackage{lmodern}
\usepackage{amssymb,amsmath,amsthm}
\usepackage{ifxetex,ifluatex}
\usepackage{fixltx2e} 

\usepackage{fullpage}
\usepackage{anysize}
\usepackage{xcolor}
\usepackage{textcomp}

\usepackage{lipsum}

\renewcommand{\S}{{\mathcal S}}
\newcommand{\A}{{\mathcal A}}
\newcommand{\R}{{\mathbb R}}

\begin{document}

 \title{Should Sports Professionals Consider Their Adversary's Strategy? A Case Study of Match Play in Golf}
    \author[aff1]{Nishad Wajge}
   \ead{nishad.wajge@gmail.com}  
    \author[aff2]{Gautier Stauffer}
   \ead{gautier.stauffer@unil.ch} 
    \address[aff1]{University of Maryland, College Park, MD 20742, United States}
      \address[aff2]{Faculty of Business and Economics (HEC Lausanne), Department of Operations, University of Lausanne, Quartier Unil-Chamberonne, 1015 Lausanne, Switzerland}

\begin{abstract}
This study explores strategic considerations in professional golf's Match Play format, challenging the conventional focus on individual performance. Leveraging PGA Tour data, we investigate the impact of factoring in an adversary's strategy. Our findings suggest that while slight strategy adjustments can be advantageous in specific scenarios, the overall benefit of considering an opponent's strategy remains modest. This confirms the common wisdom in golf, reinforcing the recommendation to adhere to optimal stroke-play strategies due to challenges in obtaining precise opponent statistics. We believe that the methodology employed here could offer valuable insights into whether opponents’ performances should also be considered in other two-player or team sports, such as tennis, darts, soccer, volleyball, etc. We hope that this research will pave the way for new avenues of study in these areas.

\end{abstract}


\maketitle

\section{Introduction}\label{intro}

Golf is a sequential game that allows players to consider their opponents’ shots and the match situation before planning their own shots. Different types of shots include the tee shot, fairway shot, approach shot, and putt. Typically played over 18 holes, golf tests a player’s decision-making, not just execution. There are two main formats of the game: Stroke Play and Match Play. Stroke Play is golf’s traditional form, where multiple golfers compete and the one with the lowest total number of shots becomes the winner. In Match Play, golfers compete head-to-head and the one winning the greatest number of holes becomes the winner. To win a hole, a golfer needs to take a smaller number of shots than its opponent on a particular hole. Therefore, the total number of shots does not really count. For a more detail description of the game of golf, we refer to \cite{GUILLOT2018}.

In Stroke Play, golfers’ decision-making strategies are mainly independent. They mainly focus on their own shots except may be for the very final holes in a tournament. In Match Play, however, golfers may be tempted to consider their opponents shots and strengths/weaknesses based on their past record before planning their own shots. It is a commonly shared piece of wisdom in the golf world that it is not a good strategy, and every player should focus on “playing his or her game”. Several professional golfers and coaches have emphasized this point over the years.  

In this paper, we challenge the corresponding wisdom through extensive numerical analysis using actual professional golfers’ data. Our goal is to analyze whether including the adversary’s statistics and position in the decision has a strong impact on the expected advantage on a hole in match-play. We focus exclusively on putting, which although involves less variations, accounts for almost half of the shots in a match. For this purpose, we develop different optimization models to retrieve the optimal strategy of a player in stroke play and in match play and we compare the performances of the former to the later. Of course, in principle, the former is suboptimal but we use the corresponding models and some simulation to quantify the corresponding gap. The optimization models are based on Markov Decision Processes and 2-Player Turned-Based Stochastic Games. Now we review some of the literature on the use of such models for strategy optimization in sports.

\subsection*{Literature review}\label{lit}

Decision-making in sports has attracted extensive research due to the large payoffs. The field is known as Sport Analytics. One of the most popular examples is celebrated in Moneyball, a book \cite{Lewis2004} (then a movie) recounting how the Oakland Athletics baseball team competed and defeated some of the wealthier national teams through the use of analytics in the early 2000’s. Statistics, data mining, machine learning, simulation, optimization, and other data science techniques are central to the field and have been used in many different sports for very different purposes, ranging from predicting scores, analyzing the effect of pressure on performance, and designing an investment strategy for the player market (see \cite{Millerbook,SportAnalyticsSI} for a recent book and a special issue on the topic; the articles published in the Journal of Sport Analytics \cite{JSA} (created in 2015) illustrate the diversity of the applications). To evaluate and enhance performance in sports, Markov chains and Markov decision processes are preferred models since they accurately represent the probability-based outcomes of each "move" executed by players or teams. These methods have been applied in various sports like tennis, basketball, volleyball, ice-hockey, golf, soccer, darts, and snooker, as seen in references e.g.~\cite{Terroba2013,trumbelj2012,Routley2015,Pfeiffer2010,Hoffmeister2015,Hoffmeister2017,Heiny2014,Maher2012,BroadieKo,Sugawara}.  

The introduction of the Shotlink\textsuperscript{\texttrademark} intelligence program has boosted academic research in golf analytics in the past 15 years, while Broadie’s  \cite{Broadie2008,Broadie2012}  strokes-gained method has revolutionized the analyses of the performance of professional golfers on the PGA Tour. In addition, a plethora of studies exploit the Shotlink\textsuperscript{\texttrademark} database to study various aspects of the game of golf (such as the effect of luck, pressure on performance, the existence of the hot hand phenomenon) through statistical analyses (e.g. \cite{Baugher2016,Ozbeklik2017,Fearing2011,Connolly2012,Robertson2014,Stockl2018,Connolly2008,Connolly2009,Connolly2012b,Hickman2015,Hickman2019,Arkes2016,Heiny2012}), performance prediction through machine learning (e.g. \cite{Hucaljuk2011,Huang2010,Moorthy2013,Wiseman2016,Lim2017,Drappi2018}, see  \cite{Bunker2019} for a recent survey), and the evaluation of different parameters (distance, dispersion, hole size) on performance through simulation and/or optimization  \cite{BansalBroadie,BroadieKo}. 
The optimization of golf strategy was first tackled by  \cite{Sugawara}. They introduced a method that blended a skill model, simulation, and Q-learning to estimate a player's optimal strategy. This is somewhat reminiscent of the model from \cite{BroadieKo}, where the assumption is that golfers always pick their theoretically best shot - a tactic frequently employed by beginners.  \cite{Sugawara} illustrated their method by applying it to "average" players, defining their competencies through parameterized distributions based on data about such players. \cite{GS2023} builds on this foundation, but with a twist. Instead of general "average" players, the study focusses on the specific empirical distribution of each PGA tour player, drawing from historical records in the Shotlink\textsuperscript{\texttrademark} database.
The research demonstrates that the inherent Markov Decision Problem linked to strategic optimization (echoing the one discussed in \cite{Sugawara}) can be solved accurately within a practical span and has a strong potential for game improvement. The corresponding methodology was designed primarily for stroke play, although it can also be used for match-play. In this research, we utilize 2-player turn-based stochastic games to identify the optimal strategy for match-play. To our knowledge, this is the first instance of these models being employed for strategy optimization in sports. One contributing factor to their limited use is the necessity for precise data regarding the opponent's gameplay, which is often challenging to acquire. This study delves into the competitive edge a player might achieve by leveraging such models.

\section{Stochastic shortest path games}

Stochastic shortest path games are 2-player extension of the Stochastic Shortest Path problem (SSP). We start with an informal description of the problem and we will give formal definitions later. The SSP is a Markov Decision Process (MDP) that generalizes the classic deterministic shortest path problem and was introduced by Bertsekas and Tsitsiklis  \cite{BertsekasTsitsiklis}. We want to control an agent, who evolves dynamically in a system composed of different {\em states}, so as to converge to a predefined  {\em target}. The agent is controlled by taking {\em actions} in each time period\footnote{We focus here on discrete time (infinite) horizon problems.}~: actions are associated with costs (possibly negative) and transitions in the system are governed by probability distributions that depend exclusively on the previous action taken and are thus independent of the past. We focus here on finite state/action spaces. The goal is to choose an action for each state, a.k.a. a {\em deterministic and stationary policy}\footnote{as for many MDPs, one can restrict to such policies, see for instance \cite{Bertsekasbook1,Bertsekasbook2,GUILLOT2018}}, so as to minimize the total expected cost incurred by the agent before reaching the (absorbing) target state, when starting from a given initial state. In order for the total expected cost to be well-defined, we need to ensure that we cannot loop in the system indefinitely while accumulating negative cost. Bertsekas and Tsitsiklis \cite{BertsekasTsitsiklis} restricted to instances where any {\em improper} policy (one that would loop indefinitely, also known as {\em transition cycle}) accumulates infinite cost. They proved that the problem can be solved using the standard techniques from MDPs, namely, {\em value iteration}, {\em policy iteration} and {\em linear programming} (see \cite{Bertsekasbook1,Bertsekasbook2,Puterman} for a detailed treatment of MDPs). Lately Bertsekas and Yu \cite{Bertsekas16} extended the assumptions to deal with improper policies of cost zero through a perturbation argument. Even more recently we proposed a slightly broader framework to deal with improper policies of cost zero by mean of polyhedral analysis \cite{GUILLOT2018}. In this project we restrict attention to SSPs where all policies are {\em proper} from any starting state, i.e. when there is no transition cycle. This setting is a special case of all previously mentioned frameworks for SSP, named SSP with inevitable termination. This is particularly convenient when dealing with the following two player extension, called {\em the Stochastic Shortest Path Game (SSPG) with inevitable termination} \cite{PatekBertsekas}\footnote{The framework proposed by Patek and Bertsekas \cite{PatekBertsekas} is more general as it allows for the existence of improper policies and simultaneous decisions of the two players (a.k.a. simultaneous stochastic games). While they prove that generalization of value iteration and policy iteration methods for MDPs converge in this setting, they do not discuss the time complexity (in particular they do not provide LP formulations).}, as both players behave `symmetrically' in this case. 

The game is played on an instance of SSP with inevitable termination but the states are now partitioned into two sets that are controlled respectively by two different players, called {\sc Min} and {\sc Max}, with antagonist objectives. The goal of {\sc Min} is to find a {\em strategy} (a choice of action for each state controlled by this player) to reach the target state with minimum expected cost (against any strategy of {\sc Max}) while {\sc Max} wants to find a strategy that maximizes the expected cost (against any strategy of {\sc Min}). This game is a special case of (zero-sum) stochastic games introduced originally by Shapley for discounted problems \cite{Shapley} but whose definition has been extended later to undiscounted problems (for a comprehensive treatment of stochastic games, see for instance \cite{NeymanSorin} and \cite{FilarVrieze}). SSPG with inevitable termination are special cases of {\em BWR-games with total effective payoff} \cite{BorosGurvich}. In particular, there exists (at least) a pair of {\em uniformly}\footnote{i.e. the policy is the same for any starting state} deterministic and stationary strategies for both players which forms a Nash Equilibrium (i.e. no player can benefit from deviating from his strategy) and the corresponding strategy for {\sc Min} minimizes the maximum expected total cost over all possible strategy for {\sc Max} and, vice-versa, the strategy for {\sc Max} maximizes the minimum expected total cost under all possible strategy for {\sc Min}. The stochastic shortest path game is the problem of finding such a pair of strategies. This framework encapsulates both (discounted) 2-player turn-based stochastic games \cite{HansenMiltersenZwick}, for which polynomial time algorithms exists when the discount factor is {\em fixed}, and {\em stopping simple stochastic games} \cite{Condon}, for which no polytime algorithm is known. 

We now give a formal definitions of the stochastic shortest path problem and the stochastic shortest path game with termination inevitable.

A stochastic shortest path instance is defined by a tuple $(\mathcal{S}, \mathcal{A},J,P,c)$ where $\S=\{0,1,\ldots,n\}$ is a finite set of {\em states},  $\A=\{0,1,\ldots,m\}$ is a finite set of {\em actions}, $J$ is a 0/1 matrix with $m$ lines and $n$ columns and general term $J(a,s)$, for all $a\in \{1,...,m\}$ and $s\in \{1,...,n\}$, with $J(a,s)=1$ if and only if action $a$ is available in state $s$, $P$ is a {\em row substochastic matrix} with $m$ lines and $n$ columns and general term $P(a,s):=p(s|a)$ (probability of ending in $s$ when taking action $a$), for all  $a\in  \{1,...,m\}$, $s\in  \{1,...,n\}$, and a cost vector $c \in \R^{m}$. The state $0$ is called the {\em target} state and the action $0$ is the unique action available in that state. Action $0$ lead to state $0$ with probability $1$. 
In the following, we denote by $\A(s)$ the set of actions available from $s\in \{1,...,n\}$ and we assume without loss of generality\footnote{If not we simply duplicate the actions.} that for all $a \in \mathcal{A}$, there exists a unique $s$ such that $a \in \mathcal{A}(s)$ (this is why we can assume that the probability of ending in $s$ when taking action $a$ is independent of the current state). We denote by 
$s(a)$ the unique state in which action $a$ is available.

A {\em stationary policy} $\Pi$ is a function that maps each state with a probability distribution over the actions. It can be represented by a $n\times m$(row) stochastic matrix $\Pi$ satisfying $\Pi(s,a)>0$ only if $J(s,a)=1$ for all $s\in \{1,...,n\}$ and $a\in \{1,...,m\}$. We will often abuse notations and use $\Pi$ to represent both a matrix and a function (this will be clear from the context). A stationary policy is {\em deterministic} if $\Pi$ is a 0/1 matrix. A stationary policy $\Pi$ is said to be {\em proper} if ${\bf 1}^T (P^T \cdot \Pi^T)^n \cdot e_i < 1$ for all $i=1,...,n$, that is, after $n$ periods of time, the probability of reaching the target state is positive, from any initial state $i$. Note that a proper stationary policy induces an absorbing Markov Chain with transition matrix $Q=P^T \cdot \Pi^T$. In particular $I-Q$ is invertible and $(I-Q)^{-1}=\lim_{K\rightarrow +\infty} \sum_{k=0}^K Q^k$ (see for instance \cite{kemeny1960finite} for more details on absorbing Markov Chains). $(Q^k e_i)(s)$ is the probability of being in state $s$ in period $k$ following policy $\Pi$ if we started in period $0$ in state $i$ and $(\Pi^T \cdot Q^k e_i)(a)$ is the probability of using action $a$ in period $k$ following policy $\Pi$ from the same starting state.  We can thus define for each $i\in \S\setminus\{0\}$, $V^\Pi(i):=\lim_{K\rightarrow +\infty}  \sum_{k=0}^{K} c^T \Pi^T Q^k e_i= c^T \Pi^T (I-Q)^{-1} e_i$ to be the expected value of policy $\Pi$ starting from state $i$. We also define $V^*(i):=\min \{V^\Pi(i) : \Pi $ proper, deterministic and stationary policy$\}$. Bertsekas and Tsistiklis \cite{BertsekasTsitsiklis} defined a stationary policy $\Pi^*$ to be {\em optimal} if $V^*(i):=V^{\Pi^*}(i)$ for all $i\in \S\setminus \{0\}$ and they proved that optimal, proper, deterministic and stationary policies exist in particular when termination is inevitable (such policies are sometimes refered to as {\em uniformly optimal} as they are optimal for any starting state). They introduced the {\em Stochastic Shortest Path Problem} as the problem of finding such an optimal deterministic and stationary policy.

An instance of a SSPG with termination inevitable is defined by a tuple $(\mathcal{S}_1,  \mathcal{S}_2 , \mathcal{A},J,P,c)$ where  $(\mathcal{S}:=\mathcal{S}_1 \cup  \mathcal{S}_2 , \mathcal{A},J,P,c)$ is an instance of SSP with inevitable termination and $\mathcal{S}_1\cap \mathcal{S}_2=\{0\}$. Note that we can define $\mathcal{A}_1=\{a\in \mathcal{A}: s(a) \in \mathcal{S}_1\}$ and $\mathcal{A}_2=\{a\in \mathcal{A}: s(a) \in \mathcal{S}_2\}$ and because we again assume w.l.o.g. that actions are available in exactly one state, $\mathcal{A}_1 \cup \mathcal{A}_2 =\A$ and $\mathcal{A}_1 \cap \mathcal{A}_2 =\{0\}$. A {\sc Min} player controls the actions in the states $\mathcal{S}_1$, while a {\sc Max} player controls the actions in the states $\mathcal{S}_2$. A (positional) {\em strategy} for player {\sc Min} is a function $\Pi_1:\S_1\mapsto \A_1$ that maps an action of $\A(s)$ to each state $s\in \S_1$,  and a (positional) {\em strategy} for player {\sc Max} is a function $\Pi_2:\S_2\mapsto \A_2$ that maps an action of $\A(s)$ to each state $s\in \S_2$. The pair $(\Pi_1,\Pi_2)$ is called a (positional) {\em strategy profile}. A strategy profile $\Pi=(\Pi_1,\Pi_2)$ induces a policy $\Pi=(\Pi_1,\Pi_2)$ for the SSP instance defined by $(\mathcal{S}, \mathcal{A},J,P,c)$ and we define the {\em value} of the strategy profile $(\Pi_1,\Pi_2)$, from an initial state $i$, as the value of the SSP solution associated with $\Pi$ i.e. $V^\Pi(i)$.

We denote by $\Sigma^1$ the set of all (positional) strategies for player {\sc Min} and by $\Sigma^2$ the set of all (positional) strategies for player {\sc Max}.  SSPG with termination inevitable are a special case of BWR-Games with total effective payoff \cite{BorosGurvich}. Because all policies are proper, for all initial state $i$, the mean payoff version\footnote{which accounts for the average cost per period over the whole horizon} of the game starting in state $i$ has value zero. It then follows from Theorem 27 in \cite{BorosGurvich} that there exists a {\em (uniform) Nash Equilibrium in positional strategy} i.e. there exists a strategy profile $(\Pi_1,\Pi_2)$ such that $V^{(\Pi_1,\Pi'_2)} (i) \leq V^{(\Pi_1,\Pi_2)} (i) \leq V^{(\Pi'_1,\Pi_2)} (i)$ for all $i\in \S\setminus \{0\}$, for all $\Pi'_1\in \Sigma^1$ and all $\Pi'_2\in \Sigma^2$.  We say that  $\Pi_1$ is the {\em best response} to strategy $\Pi_2$ and $\Pi_2$ is the {\em best response} to strategy $\Pi_1$. Now by Von Neumann's minimax theorem for zero-sum games \cite{VonNeumann}, we know that such a  Nash Equilibrium  $(\Pi_1,\Pi_2)$ satisfies $$\displaystyle V^{(\Pi_1,\Pi_2)} (i) = \min_{\Pi'_1\in \Sigma^1} \max_{\Pi'_2\in \Sigma^2} V^{(\Pi'_1,\Pi'_2)}(i)= \max_{\Pi'_2\in \Sigma^2} \min_{\Pi'_1\in \Sigma^1} V^{(\Pi'_1,\Pi'_2)}(i), \mbox{ for all } i\in \S\setminus \{0\} $$ and moreover  $V^{(\Pi_1,\Pi_2)} (i)$ is identical for each Nash Equilibrium. $V^{(\Pi_1,\Pi_2)} (i)$ is then called the {\em value of the game} when starting from state $i$. The stochastic shortest path game (with inevitable termination) is the problem of finding such a Nash equilibrium (or the value of the game). There is no known polynomial time algorithm for this problem but strategy iteration is a combinatorial algorithm that solves the problem efficiently in practice (and is guaranteed to converge in a finite number of steps).


\begin{algorithm}
\caption{Strategy Iteration}
\begin{algorithmic}[1]
\Procedure{StrategyIteration}{}
    \State Set \(\Pi(s)\) to a random action from $\mathcal{A}$ (available in $s$) for each state \(s\) in $\mathcal{S}$
    \While{True}
        \State  \text{Policy Evaluation} : compute $V^{\Pi}$ 
        \State {\bf Optimize for {\sc Max} while fixing {\sc Min}'s policy:}
        \While{$\exists$  \(s\) in \(\mathcal{S}_2\) and  $\exists$ action \(a \in \mathcal{A}_2(s) \) : $c_a + \sum_{s'\in \mathcal{S}} P(s' | a ) V^{\Pi}(s') > V^{\Pi}(s)$}
               \State \(\Pi(s) = a\)
               \State  \text{Policy Evaluation} : compute $V^{\Pi}$ 
        \EndWhile    
        \State {\bf Check improvement for {\sc Min}:}
          \If{$\exists$  \(s\) in \(\mathcal{S}_1\) and  $\exists$ action \(a \in \mathcal{A}_1(s) \) : $c_a + \sum_{s'\in \mathcal{S}} P(s' | a ) V^{\Pi}(s') < V^{\Pi}(s)$}
                \State \(\Pi(s) = a\)
          \Else
                \State \text{Break}
          \EndIf
    \EndWhile
\EndProcedure
\end{algorithmic}
\end{algorithm}

\section{Modelling putting}

In this study, we primarily examine the art of putting in golf. As noted previously, putting usually constitutes about half of a golfer's shots. To streamline our analysis, we only consider greens that are flat. Our primary objective is to minimize the expected number of putts a golfer takes once on the green. The crucial factor in this context is the distance to the hole, rather than the specific location on the green.

For the subsequent analysis, we adopt a basic skill model for putting. A player selects a target direction and distance, acknowledging the possibility of errors in either choice. Suppose a player is positioned at $O = (0,0)$ and aims for the destination $T = (0,d)$. We postulate that the ball's final position (assuming no obstacles) is a random variable, $X_d$, which follows a 2D probability density function. To elaborate, we consider the sample space $\Omega := \{(x,y) \in \mathbb{R}^2\}$. We hypothesize that $X_d$ adheres to a probability density function $f_d : (x,y) \in \Omega \mapsto f_d(x,y) \in \mathbb{R}_+$ for every $d$ within the range $[0,D]$, where $D$ represents the greatest distance on a green.

\begin{figure}[h!]
\begin{tabular}{cc}
\includegraphics[width=0.8\textwidth]{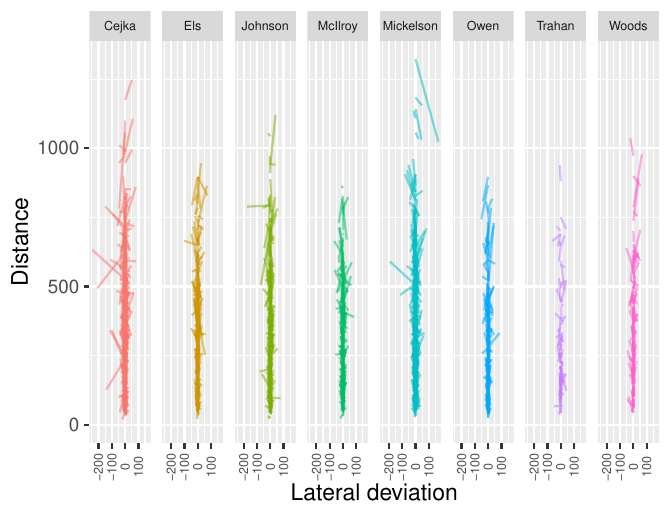} 
\end{tabular}
\caption{Putting data for a selection of PGA Tour players. }\label{fig_putting}
\end{figure}

Figure \ref{fig_putting} depicts putting data for a selected set of professional golfer. Each segment signifies the hole's position (rotated so that the hole is on the y-axis) and the putt's outcome. This distribution does not precisely represent an empirical distribution sampled from $f_d$ for several reasons. Firstly, these putts were executed on undulating greens, and some were successful (in these instances, the segment appears as a dot). Secondly, golfers rarely aim directly for the hole. A popular saying in golf, "never up, never in," suggests golfers often aim slightly beyond the hole. Nevertheless we can use these data to estimate $X_d$ by presuming the putts were made on flat greens and making a few more assumptions:

\begin{enumerate}
    \item We posit that the angle dispersion is independent of the distance and follows a normal distribution centered in zero and we estimate the standard deviation from the raw data above by considering all putts made. The value are provided in Table \ref{tab0}.
\item We estimate the target distance for a selection of hole distances by calculating the mean distance covered by the ball for all putts that either (a) missed the hole but went beyond it or (b) missed the hole and were not directly aligned with it. This approach aims to mitigate the effect of sinking a putt. Additionally, we estimate the standard deviation of the distance error for these corresponding distances. The standard deviations for these distances are detailed in Table \ref{tab1}. It is important to note that for each distance $d$, our estimates were constructed from a set of 100 putts where the initial distance to the hole was as close as possible to $d$. If the actual distance was $d' > 0$, we rescaled the final putt destination $(x, y)$ to $\frac{d}{d'}\cdot (x, y)$.

    \item We posit that the distance error follows a normal distribution centered on the targeted distance with a standard deviation that is interpolated from the value estimated in Table \ref{tab1}.
\end{enumerate}

\begin{table}[ht]
\centering
\begin{tabular}{rlrr}
  \hline
 & last name &  angle sd \\ 
  \hline
1 & Cejka  & 0.028 \\ 
  2 & Els  & 0.025 \\ 
  3 & Johnson  & 0.027 \\ 
  4 & McIlroy & 0.029 \\ 
  5 & Mickelson  & 0.028 \\ 
  6 & Owen  & 0.029 \\ 
  7 & Trahan  & 0.031 \\ 
  8 & Woods  & 0.025 \\ 
   \hline
\end{tabular}
 \caption{}\label{tab0}
\end{table}

\begin{table}[ht]
\centering
\begin{tabular}{lrrr}
  \hline
last name & hole (in.) & target (in.)  & sd (in.) \\ 
  \hline
Johnson & 40.00 & 58.33 & 15.65 \\ 
  Johnson & 100.00 & 117.84 & 12.41 \\ 
  Johnson & 200.00 & 211.65 & 17.10 \\ 
  Johnson & 400.00 & 406.10 & 33.84 \\ 
  Johnson & 800.00 & 800.70 & 55.82 \\ 
  Els & 40.00 & 58.71 & 15.15 \\ 
  Els & 100.00 & 113.76 & 9.56 \\ 
  Els & 200.00 & 212.41 & 19.29 \\ 
  Els & 400.00 & 409.76 & 33.50 \\ 
  Els & 800.00 & 802.37 & 44.59 \\ 
  Woods & 40.00 & 64.49 & 22.49 \\ 
  Woods & 100.00 & 114.81 & 11.83 \\ 
  Woods & 200.00 & 220.71 & 18.72 \\ 
  Woods & 400.00 & 412.20 & 30.38 \\ 
  Woods & 800.00 & 808.77 & 44.61 \\ 
  Mickelson & 40.00 & 67.10 & 19.01 \\ 
  Mickelson & 100.00 & 121.29 & 13.19 \\ 
  Mickelson & 200.00 & 215.84 & 18.32 \\ 
  Mickelson & 400.00 & 407.37 & 37.79 \\ 
  Mickelson & 800.00 & 801.95 & 59.84 \\ 
  McIlroy & 40.00 & 56.89 & 11.27 \\ 
  McIlroy & 100.00 & 122.02 & 15.49 \\ 
  McIlroy & 200.00 & 217.90 & 14.40 \\ 
  McIlroy & 400.00 & 411.13 & 35.01 \\ 
  McIlroy & 800.00 & 798.36 & 47.98 \\ 
  Cejka & 40.00 & 60.36 & 16.28 \\ 
  Cejka & 100.00 & 113.59 & 13.57 \\ 
  Cejka & 200.00 & 211.14 & 18.02 \\ 
  Cejka & 400.00 & 398.35 & 37.14 \\ 
  Cejka & 800.00 & 795.15 & 69.29 \\ 
  Owen & 40.00 & 63.69 & 13.17 \\ 
  Owen & 100.00 & 117.67 & 13.03 \\ 
  Owen & 200.00 & 212.63 & 15.79 \\ 
  Owen & 400.00 & 402.33 & 40.88 \\ 
  Owen & 800.00 & 798.66 & 49.42 \\ 
  Trahan & 40.00 & 54.59 & 19.68 \\ 
  Trahan & 100.00 & 120.60 & 13.02 \\ 
  Trahan & 200.00 & 213.92 & 19.95 \\ 
  Trahan & 400.00 & 408.32 & 31.98 \\ 
  Trahan & 800.00 & 808.22 & 32.72 \\ 
   \hline
\end{tabular}
\caption{Estimated targeted distance for some PGA Tour players. Hole represents the distance to the hole, target the estimated distance targeted and sd the standard deviation.}\label{tab1}
\end{table}

\begin{figure}[h!]
\centering
\includegraphics[width=0.8\textwidth]{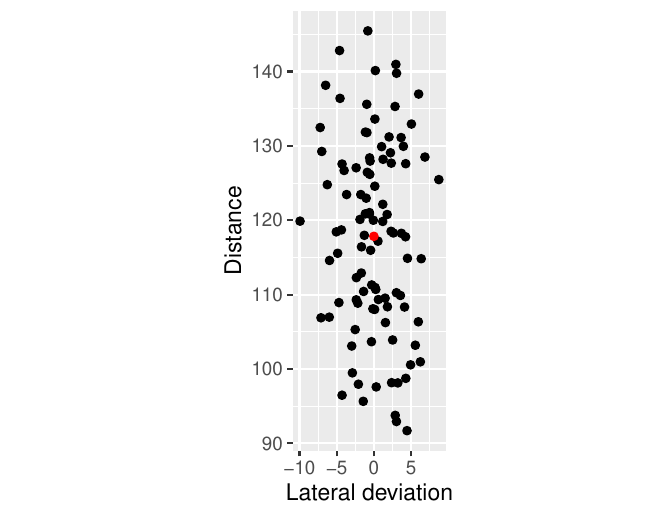}
\caption{Estimated putt distribution for a targeted distance of 117.83 inches (highlighted in red).}\label{fig2}
\end{figure}

Now, we can construct a representative empirical distribution of \(X_{d}\) for any player and any target distance \(d\). This assumes that the directional error and the distance error are independent of each other. Figure \ref{fig2} provides an illustrative example for Dustin Johnson with a targeted distance of 117.83 inches.

\begin{figure}[h!]
\begin{tabular}{cc}
\includegraphics[width=0.8\textwidth]{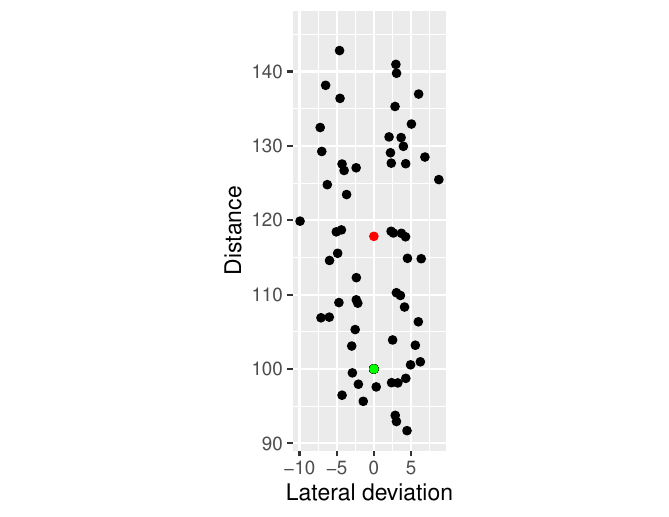} 
\end{tabular}
\caption{Estimated distribution of the final destination for a targeted distance of 117.83 inches (in red) when the hole is at a distance of 100 inches (in green). Note that in comparison to Fig. \ref{fig2}, some putts are now resting in the hole.} \label{fig3}
\end{figure}

Given the estimated putt distribution, we can simulate the putts of Dustin Johnson that would be captured by the hole, assuming the hole is placed at a distance of say 100 inches. We utilize a formula derived by Penner \cite{Penner} to make this assessment. Penner posits that a ball, approaching the hole at speed \(s\) (in m/s) and at a distance \(\delta \leq R\) from the hole's center (in m),  is effectively captured if \(s < 1.63 - 1.63 (\frac{\delta}{R})^2\), where \(R = 0.054\) m denotes the hole's radius. On a flat green, determining the lateral deviation \(\delta\) at the hole is straightforward. What remains is the calculation of the speed \(s\). The relationship between distance and speed on a flat green can be described as \(D = k \cdot s^2\), where the constant \(k\) varies based on the green's speed. For greens clocking a speed of 12 feet (a measurement derived using a stimpmeter in golf and representative of the average green speed on the PGA Tour), \(k\) equals 1.093. This relationship enables us to deduce the ball's speed at the hole based on its final resting position in the absence of obstacles. If the hole is \(d\) units away and the ball's final position is \(t\geq d\) units, then the speed at the hole is approximated as \(s = \sqrt{(t-d)/1.093}\). Implementing this model yields the results shown in Figure \ref{fig3}.
  
We simulated outcomes for our chosen group of players across various distances, based on our estimates of their target distances. These findings can be viewed in Table \ref{fig4}. While the figures may not align precisely in every case, the models do appear to provide reasonable approximations of the respective players. It is important to note that our objective in this study is not to create flawless digital replicas of the players. Instead, our goal is to generate data representative of typical PGA Tour players, and this intent is affirmed by the information in the table.

\begin{table}[ht]
\centering
\begin{tabular}{rrllrrrrr}
  \hline
  last name & CR (model)  & CR (data) & RD (model) & RD (data) & DH \\ 
  \hline
Johnson & 0.41 & 0.41 & 15.72 & 13.71 & 100.00 \\ 
  Johnson & 0.25 & 0.21 & 17.17 & 18.79 & 200.00 \\ 
  Johnson & 0.05 & 0.05 & 30.82 & 32.16 & 400.00 \\ 
  Johnson & 0.01 & 0.03 & 48.42 & 42.08 & 800.00 \\ 
  Els & 0.53 & 0.41 & 15.53 & 13.71 & 100.00 \\ 
  Els & 0.18 & 0.21 & 18.68 & 18.79 & 200.00 \\ 
  Els & 0.12 & 0.05 & 28.32 & 32.16 & 400.00 \\ 
  Els & 0.02 & 0.03 & 43.58 & 42.08 & 800.00 \\ 
  Woods & 0.46 & 0.41 & 16.46 & 13.71 & 100.00 \\ 
  Woods & 0.16 & 0.21 & 23.73 & 18.79 & 200.00 \\ 
  Woods & 0.05 & 0.05 & 31.32 & 32.16 & 400.00 \\ 
  Woods & 0.02 & 0.03 & 42.82 & 42.08 & 800.00 \\ 
  Mickelson & 0.42 & 0.41 & 25.62 & 13.71 & 100.00 \\ 
  Mickelson & 0.22 & 0.21 & 19.87 & 18.79 & 200.00 \\ 
  Mickelson & 0.02 & 0.05 & 32.59 & 32.16 & 400.00 \\ 
  Mickelson & 0.01 & 0.03 & 53.03 & 42.08 & 800.00 \\ 
  McIlroy & 0.34 & 0.41 & 25.20 & 13.71 & 100.00 \\ 
  McIlroy & 0.18 & 0.21 & 19.15 & 18.79 & 200.00 \\ 
  McIlroy & 0.09 & 0.05 & 35.82 & 32.16 & 400.00 \\ 
  McIlroy & 0.00 & 0.03 & 45.46 & 42.08 & 800.00 \\ 
  Cejka & 0.37 & 0.41 & 16.79 & 13.71 & 100.00 \\ 
  Cejka & 0.15 & 0.21 & 19.60 & 18.79 & 200.00 \\ 
  Cejka & 0.09 & 0.05 & 31.72 & 32.16 & 400.00 \\ 
  Cejka & 0.01 & 0.03 & 59.41 & 42.08 & 800.00 \\ 
  Owen & 0.33 & 0.41 & 18.24 & 13.71 & 100.00 \\ 
  Owen & 0.18 & 0.21 & 17.50 & 18.79 & 200.00 \\ 
  Owen & 0.08 & 0.05 & 34.55 & 32.16 & 400.00 \\ 
  Owen & 0.02 & 0.03 & 43.97 & 42.08 & 800.00 \\ 
  Trahan & 0.42 & 0.41 & 22.31 & 13.71 & 100.00 \\ 
  Trahan & 0.16 & 0.21 & 21.21 & 18.79 & 200.00 \\ 
  Trahan & 0.04 & 0.05 & 30.99 & 32.16 & 400.00 \\ 
  Trahan & 0.01 & 0.03 & 36.05 & 42.08 & 800.00 \\ 
   \hline
\end{tabular}
\caption{Comparison of actual and simulated data for a selection of distances and players. CR represents the percentage of ball captures (capture reate), RD represents the remaining distance to the hole and DH the distance to the hole.}\label{fig4}
\end{table}

\section{Stochastic Shortest Path models and results}

We are now ready to build strategy optimization models for the different players. Let us focus first on the optimization of the expected number of putts of a player from any position on the green. Because the green is flat, we can represent the situation of a player solely by the distance to the pin. We consider all possible positions on the green with a maximum distance of 800 inches and with a precision of $\delta=5$ inches. We thus consider a state space $\mathcal S=\{0,1,...,n\}$ for $n=160$ where state $i$ represent a putt at distance $i\cdot \delta$ (state $0$ representing the hole). Then, because we assume a player makes no aiming error, an action is essentially a selection of a target beyond the hole (there is no point choosing a target before on a flat green). 
The maximum meaningful distance beyond the hole is the distance corresponding to the maximum speed that allows to hole the putt (1.63 m/s). This corresponds to $\sim110$ inches at the speed considered above. We also consider a similar precision of $\delta=5$ inches for the action choice. So we consider an action space ${\mathcal A}=\{(s,j)\}$ for all $s\in S\setminus \{0\}$ and for $j\in \{0,1,...,m\}$ with $m=22$ (there is an additional action $0$ in the target state $0$ that lead to $0$ with probability one). Note that if a putt ends up at a distance over 800 inches, we assume it ends at a distance of 800 inches exactly (it will happen only with an extremely low probability given the fact that we target a distance at most 110 inches away from the pin and given the distance error - see Table \ref{tab1}).

In this problem, $J(s,(s',a))=1$ if and only if  $s=s' \in \{1,...,n\}$ and $a\in \{0,1,...,m\}$ and $\Pi(s',(s,a))$ is the probability of ending up at distance $s' \cdot \delta$ when targeting distance $(s+a) \cdot \delta $. This can easily be computed from the estimates of the final destinations computed as in Fig. \ref{fig3}. We have simulated 1000 putts to obtain the estimates of $\Pi(s',(s,a))$. We have assumed, as discussed already, normal distribution centered in 0 for the angle and the distance dispersions, with constant standard deviation for the angle deviation and linear interpolation for the standard deviation on the distance, taken from Tables \ref{tab0} and \ref{tab1} respectively.

\begin{figure}[h!]
\begin{tabular}{cc}
\includegraphics[width=0.8\textwidth]{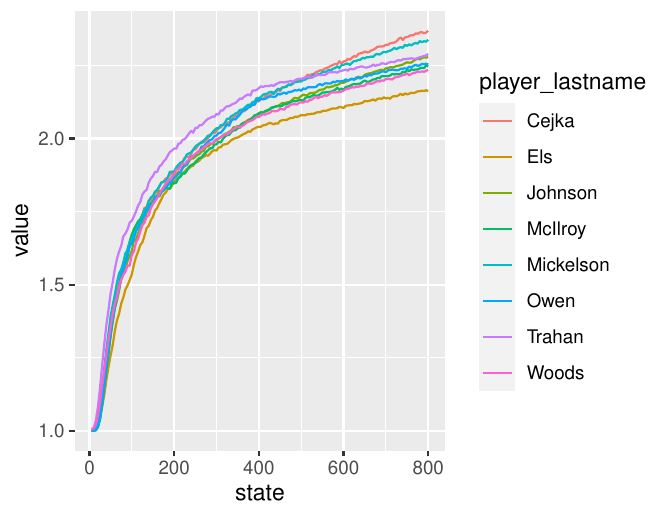} 
\end{tabular}
\caption{Estimated putting performances using the optimal stroke-play strategy}\label{fig5}
\end{figure}

We used value iteration to optimize the policy and the value function of each player is given in Fig \ref{fig5} (the reward $c((s,a))$ is one for each action except for). Again, we are not concerned here with having a perfect fit with the true performances of the corresponding players. What matters to us is that the corresponding profiles are representative of different professional performances. The figures are aligned with average performances from the literature (see Fig. 1 in \cite{James2008}) and there are interesting variations.

We can observe the optimal distance targeted beyond the hole for the various players in Fig \ref{fig6}. The figure shows that the (typical) targeted distance beyond the hole decreases with the distance, which is to be expected as players try to arbitrate between maximizing the chances of holing the putt and minimizing the expected remaining distance. Now the (apparently) hectic fluctuations in Fig \ref{fig6} mainly come from the fact that the value for nearby actions is pretty similar and distance discretization induces rounding errors. 

\begin{figure}[h!]
\begin{tabular}{cc}
\includegraphics[width=0.8\textwidth]{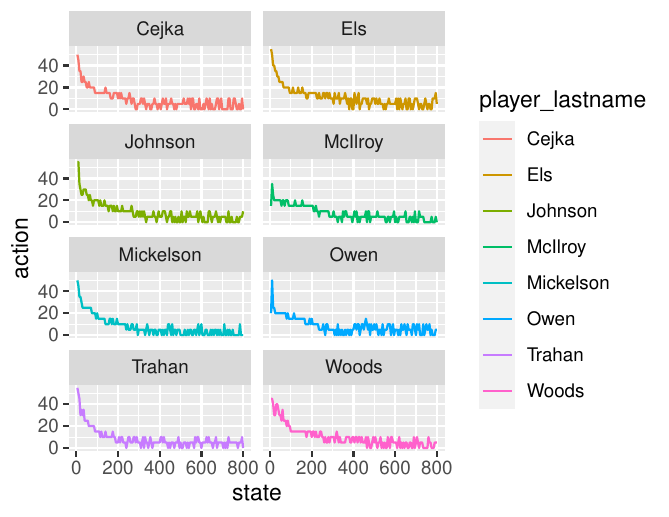} 
\end{tabular}
\caption{Estimated putting performances}\label{fig6}
\end{figure}

We are now ready to study the match play setting. Consider two players names player 1 and player 2 respectively. We assume that both players want to maximize the expected number of points they can get in any situation on a green. Now a situation is characterized by the position of each player and the number of point difference. Note that we can limit the number of point difference to $5$ as whatever the position of the two players on the green, it is hardly possible (especially at the professional level) for a player to win or tie if the adversary is 5 shot ahead (the most favorable situation would be that the player is in the hole already and the other would make 5 putts but the probability that a PGA tour player makes 4 putts or more is very low : 3 cases out of 10000 putts in our data set). We can define an instance of SSPG with termination inevitable as follows. 

The different states $\mathcal{S}$ in the system are all triplets $(s,s',\Delta)$ for $s,s'\in \{0,1,...,n\}$ and $\Delta \in \{-5,...,5\}$: $s$ represent the position of player 1, i.e. ball at distance $s\cdot \delta$, $s'$ the position of player 2, i.e. ball at distance $s'\cdot \delta$, and $\Delta$ represents the difference in number of shots between player 1 and player 2. We assign each triplet $(s,s',\Delta)$ to $\mathcal{S}_1$ or $\mathcal{S}_2$ following the rule of golf (the player further away plays first): $(s,s',\Delta) \in \mathcal{S}_1$ if $s>s'$, $(s,s',\Delta) \in \mathcal{S}_2$ if $s<s'$ and $(s,s',\Delta)$ is assigned randomly if $s=s'$. We slightly abuse our definition here to simplify the exposition and we consider several terminal states in this model: each state  $(s,s',\Delta)$ where $\Delta=5$ is a terminal state with value $-1$ (player 1 loses) and each state  $(s,s',\Delta)$ where $\Delta=-5$ is a terminal state with value $1$ (player 2 loses). It is in fact similar to having a unique action in each of these states that would lead to target state $0$ with the reward equal to the value. From a state $(s,s',\Delta) \in \mathcal{S}_1$, we consider all actions of player 1, that is, all pairs $((s,s',\Delta),j)$ for $j\in \{0,1,...,m\}$. When taking action $((s,s',\Delta),j)$, we end in state $(s'',s',\Delta+1)$ with probability $p$ where $p$ is the probability that player 1 ends in state $s''$ when playing action $j$ in his or her stroke play SSP model. Similarly, from a state $(s,s',\Delta) \in \mathcal{S}_2$, we consider all actions of player 2, that is, all pairs $((s,s',\Delta),j)$ for $j\in \{0,1,...,m\}$. When taking action $((s,s',\Delta),j)$, we end in state $(s,s'',\Delta-1)$ with probability $p$ where $p$ is the probability that player 2 ends in state $s''$ when playing action $j$ in his or her stroke play SSP model. Player 1 is the max player and player 2 is the Min player here.

The following figures illustrate the difference between stroke play and optimal match-play strategy for the second player across a (random) set of 9 match-play scenarios involving the 8 players above (with the first player employing the optimal match-play strategy). In these figures, red represents a more aggressive optimal policy (targeting a distance at least 10 inches further), while green signifies the opposite. We plot the differences for different shot (dis-)advantage on the green. If a plot is not given, it means that there is no difference.

\begin{figure}[h!]
\begin{tabular}{cc}
\includegraphics[width=0.8\textwidth]{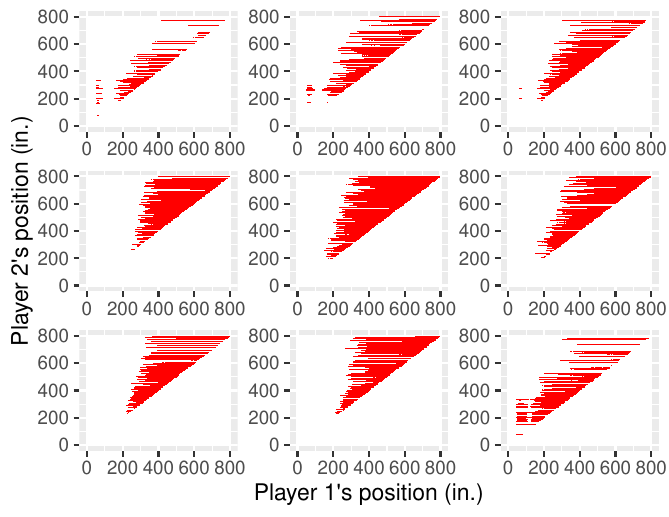} 
\end{tabular}
\caption{Difference in strategy when player 2 has a 2-shot disadvantage}\label{fig6}
\end{figure}

\begin{figure}[h!]
\begin{tabular}{cc}
\includegraphics[width=0.8\textwidth]{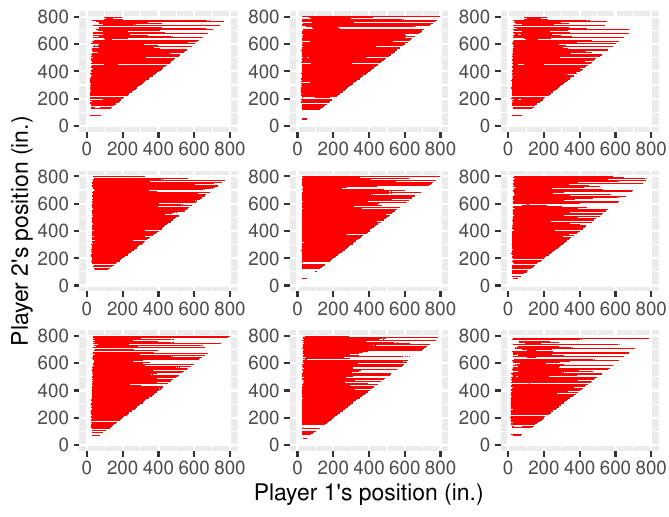} 
\end{tabular}
\caption{Difference in strategy when player 2 has a 1-shot disadvantage}\label{fig6}
\end{figure}

\begin{figure}[h!]
\begin{tabular}{cc}
\includegraphics[width=0.8\textwidth]{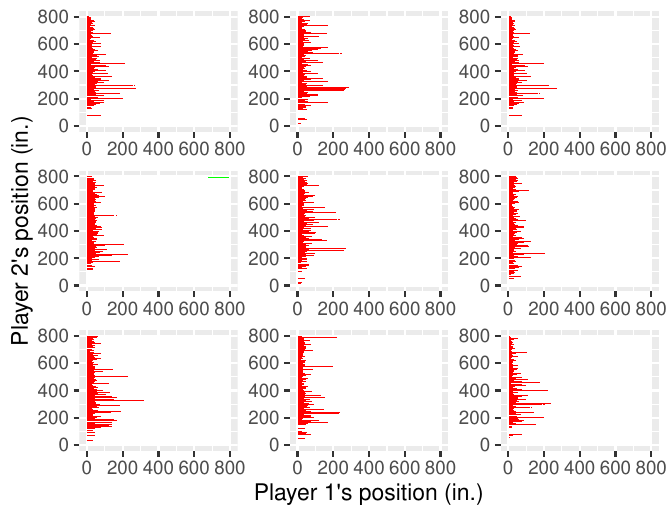} 
\end{tabular}
\caption{Difference in strategy when both players reach the green in the same number of shots.}\label{fig6}
\end{figure}

\begin{figure}[h!]
\begin{tabular}{cc}
\includegraphics[width=0.8\textwidth]{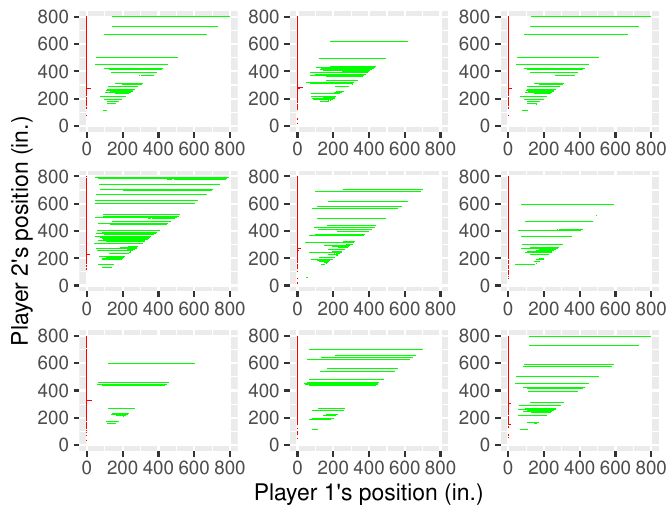} 
\end{tabular}
\caption{Difference in strategy when player 2 has a 1-shot advantage}\label{fig6}
\end{figure}

\begin{figure}[h!]
\begin{tabular}{cc}
\includegraphics[width=0.8\textwidth]{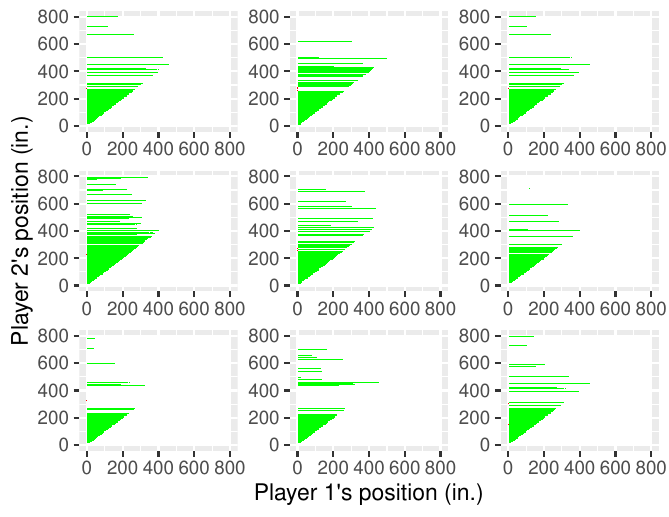} 
\end{tabular}
\caption{Difference in strategy when player 2 has a 2-shot advantage}\label{fig6}
\end{figure}

\begin{figure}[h!]
\begin{tabular}{cc}
\includegraphics[width=0.8\textwidth]{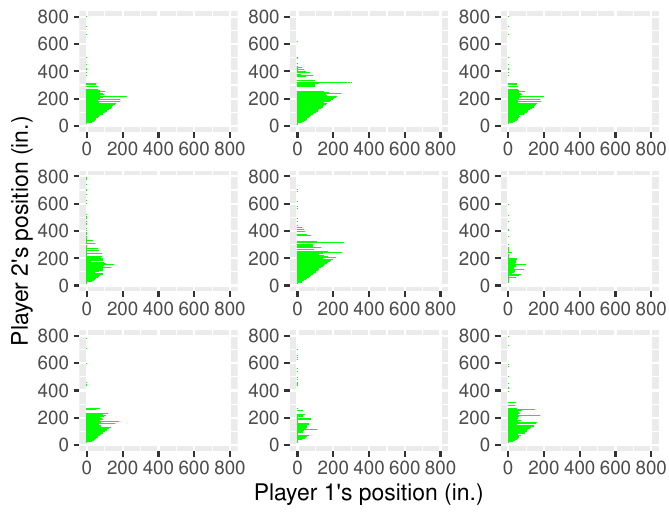} 
\end{tabular}
\caption{Difference in strategy when player 2 has a 3-shot advantage}\label{fig6}
\end{figure}

The figures clearly illustrate significant differences between stroke-play and match-play strategies. It is worth noting that the occasional, non-smooth "stripes" are a result of two nearby locations having sometimes up to a 10-inch variation in the optimal stroke-play strategy. This discrepancy arises because the optimal stroke-play values are nearly identical for nearby positions and we face some rounding effects as already discussed in Fig. \ref{fig6}. Nevertheless, the trends conveyed by the figures are quite evident.

In situations where a player is at a disadvantage (except when the disadvantage is insurmountable), a slightly more aggressive approach is advisable. Conversely, when a player is in an advantageous position, a somewhat more conservative strategy is preferable. The most intriguing insight emerges when players have reached the green in the same number of shots. In such cases, there exists a threshold for player one's position. Below this threshold, player two should adopt a more aggressive approach, attempting to compensate for the distance disadvantage with a higher probability of sinking the putt.

Turning our attention to the objective function, the results are somewhat mixed. Table \ref{tabU} presents the mean and maximum expected benefit for player two across the nine matches considered earlier. If player two were fortunate enough to consistently achieve the best possible benefit, they would gain an additional 0.8028 points over a golf course with an average benefit falling within the range of $[0.0486, 0.117]$, depending on how often player two finds themselves at an advantage, disadvantage, or neutral position in terms of shot difference.

These numbers are relatively modest and are unlikely to provide a significant competitive advantage to player two. Given the challenges in collecting precise statistics regarding opponents, we would advise professional players to adhere to their optimal stroke-play strategy. Nevertheless, we recommend the continued collection of accurate data regarding their game, which can be used to determine their stroke-play strategy, particularly in putting (following the methodology outlined here) or for their overall game, as discussed in \cite{GUILLOT2018}.

\begin{table}[ht]
\centering
\begin{tabular}{rrr}
  \hline
point difference & mean value & max value \\ 
  \hline
-5 & 0.0000 & 0.0000 \\ 
  -4 & 0.0000 & 0.0000 \\ 
  -3 & 0.0000 & 0.0000 \\ 
  -2 & 0.0015 & 0.0141 \\ 
  -1 & 0.0065 & 0.0355 \\ 
  0 & 0.0027 & 0.0446 \\ 
  1 & 0.0063 & 0.0446 \\ 
  2 & 0.0084 & 0.0434 \\ 
  3 & 0.0029 & 0.0435 \\ 
  4 & 0.0003 & 0.0228 \\ 
  5 & 0.0000 & 0.0000 \\ 
   \hline
\end{tabular}
\caption{The table represents the mean and max gap between the average optimal match-play strategy and the stroke-play strategy for different point difference between the players (e.g. a point difference of -1 indicates that player 1 is on the green with one shot less than player 2. The mean value is averaged over all players and distances. The mean and max values are similar for each individual players.}\label{tabU}
\end{table}

\section{Conclusion and perspective}

Our preliminary analysis indicates that there is only a marginal overall advantage in employing the optimal match-play strategy over the stroke-play strategy, particularly when it comes to putting (excluding obvious situations). This finding aligns with the conventional wisdom in golf, especially concerning this specific aspect of the game. We anticipate that similar conclusions would likely emerge for other aspects of the game, although further investigations would be necessary to confirm this as other types of golf shots, such as tee shots and approach shots, inherently exhibit greater variability than putting. 

In an interview, major champion Collin Morikawa attributed his early success to what he calls a "mastery mindset." This mindset empowers him to maintain unwavering focus on his own game, without becoming distracted by the performance of others. Our study underscores the merit of this approach, especially given the additional potential psychological impact of constantly reacting and adjusting one's strategy to the opponent's game.

We believe that the methodology employed here could offer valuable insights into whether opponents' performances should also be considered in other two-player or team sports, such as tennis, darts, soccer, volleyball, etc. We hope that this research will pave the way for new avenues of study in these areas.

\section{Acknowledgement}

Nishad Wajge would like to thank Akshat Wajge for helping debug complex R and Python programs and Rajesh Wajge for helpful brainstorming on dynamic programming and golf. 
We would like to thank the PGA Tour for giving us access to the ShotLink\textsuperscript{\texttrademark} data.

\bibliographystyle{plain} 
\bibliography{biblio_golf}  

{

\end{document}